\documentclass[12pt,a4paper]{article}
\usepackage[colorlinks=true,linkcolor=black, citecolor=blue, urlcolor=blue]{hyperref}
\usepackage{epsfig}
\usepackage{graphics}
\usepackage{amsfonts}
\usepackage{amsmath}
\usepackage{latexsym,mathrsfs}

\newtheorem{theorem}{Theorem}[section]

\newtheorem{proposition}{Proposition}[section]

\newtheorem{remark}{Remark}[section]

\DeclareMathOperator{\supp}{supp} 
 \numberwithin{equation}{section}

\title{\textbf{On the projections of the multifractal Hewitt-Stromberg dimension functions}}
\author{Bilel Selmi}

\begin{document}
\maketitle
\begin{abstract}
The aim of this paper is to study the behavior of the multifractal
Hewitt-Stromberg dimension functions under projections in Euclidean
space. As an application, we study the multifractal analysis of the
projections of a measure. In particular,  we obtain general results
for the multifractal analysis of the orthogonal projections on
$m$-dimensional linear subspaces of a measure $\mu$ satisfying the
multifractal formalism which is based on the Hewitt-Stromberg
measures.
\end{abstract}

\noindent{\bf MSC-2010:} {28A20, 28A78, 28A80}


\section{Introduction and statement of the results}

Recently, the projection behavior of dimensions and multifractal
spectra of sets and measures have generated large interest in the
mathematical literature \cite{BF, DB, FM1, FJ, JJ, J, J1, SP, SP1}.
The multifractal analysis is always interesting therefore studies
are developing: especially the case of projections of measures, as
we did not know till now an exact link between the dimension of the
original object (measure, set,...) and its projections: How to
relate the dimension of the measure with its projections? It is for
example still not known by projecting a multifractal measure on a
multifractal set, such as Cantor's types, how to deduce the
dimension of the projection measure from the properties of the
original one and the dimensions of the projection and the set we
project on it. Dimensional properties of projections of sets and
measures have been investigated for decades. The first significant
work in this area was the result of Marstrand \cite{M}, to which the
Hausdorff dimension of a planar set is preserved under typical
orthogonal projections. This result was later generalized to higher
dimensions by Kaufman \cite{K}, Mattila \cite{MP}, and Hu and Taylor
\cite{HT}. They obtained similar results for the Hausdorff dimension
of a measure. It is natural to seek projection results for other
notions of dimension. However, examples show that the direct
analogue of Marstrand's projection theorem is not valid for
box-counting and packing dimensions, though there are non-trivial
lower bounds on the dimensions of the projections. The fact that the
box-counting and packing dimensions of the projections of a Borel
set $E$ are constant for almost all subspaces $V$ was established by
Falconer and Mattila in \cite{FM} and Falconer and Howroyd in
\cite{Falconer1, FH, FH1, JDH}. Such constant value, given by a
dimension profile of $E$, was specified somewhat indirectly. For
packing dimensions, this is given in terms of the suprema of
dimension profiles of measures supported by $E$ which in turn are
given by critical parameters for certain almost sure pointwise
limits \cite{FH}. The approach in \cite{FH1, JDH} defines
box-counting dimension profiles in terms of weighted packings
subject to constraints. However, despite these substantial advances
for fractal sets, only very little is known about the multifractal
structure of projections of measures, except a paper by O'Neil
\cite{O}. Later, in \cite{BB} Barral and Bhouri, studied the
multifractal analysis of the orthogonal projections on
$m$-dimensional linear subspaces of singular measures on
$\mathbb{R}^n$ satisfying the classical multifractal formalism. The
result of O'Neil was later generalized by Selmi et al. in \cite{DB,
BD1, BD3, SBB1, SB2}.

The notion of singularity exponents or spectrum and generalized
dimensions are the major components of the multifractal analysis.
They were introduced with a view of characterizing the geometry of
measure and to be linked with the multifractal spectrum which is the
map which affects the Hausdorff or packing dimension of the
iso-H\"{o}lder sets
$$
E_\mu(\alpha)=\left\{x\in \supp\mu\;\Big|\; \lim_{r\to
0}\frac{\log\big(\mu B(x,r)\big)}{\log r}=\alpha \right\}
$$
for a given $\alpha\geq 0$ and $\supp\mu$ is the topological support
of probability measure $\mu$ on $\mathbb{R}^n$, $B(x,r)$ is the
closed ball of center $x$ and radius $r$. It unifies the
multifractal spectra to   the multifractal packing function
$b_\mu(q)$ and the multifractal packing function $B_\mu(q)$ via the
Legendre transform \cite{NB2, NB1, BBH, OL}, i.e.,
$$
\dim_H\Big(E_\mu(\alpha)\Big)=\inf_{q\in\mathbb{R}}\Big\{q\alpha+b_{\mu}(q)\Big\}=:
b_{\mu}^*(\alpha)
$$
and
$$
\dim_P\Big(E_\mu(\alpha)\Big)=\inf_{q\in\mathbb{R}}\Big\{q\alpha+B_{\mu}(q)\Big\}=:
B_{\mu}^*(\alpha).
$$
There has been a great interest in understanding the fractal
dimensions of projections of the iso-H\"{o}lder sets and measures.
The authors in \cite{O, SB3}  compared the multifractal packing
function $b_\mu(q)$ and the multifractal packing function $B_\mu(q)$
of a set $E$ of $\mathbb{R}^n$ with respect to a measure $\mu$ with
those of their projections onto $m$-dimensional subspaces. And they
obtained general results for the multifractal analysis of the
orthogonal projections on $m$-dimensional linear subspaces of a
measure $\mu$ satisfying the classical multifractal formalism which
is based on the Hausdorff and packing measures on
$\gamma_{n,m}$-almost every linear $m$-dimensional subspaces for
$q>1$, where $\gamma_{n,m}$ is the uniform measure on $G_{n,m}$, the
set of linear $m$-dimensional subspaces of $\mathbb{R}^n$ endowed
with its natural structure of a compact metric space, see \cite{MP1}
for precise definitions of this.

Hewitt-Stromberg measures were introduced  in \cite[Exercise
(10.51)]{HeSt}. Since then, they  have  been investigated  by
several authors, highlighting their importance in the study of local
properties of fractals and products of fractals. One can cite, for
example \cite{Ha1, Ha2, JuMaMiOlSt, Olll2, Zi}. In particular,
Edgar's textbook \cite[pp. 32-36]{Ed} provides an excellent and
systematic introduction to these  measures. Such measures appear
also explicitly, for example, in  Pesin's monograph \cite[5.3]{Pe}
and implicitly in Mattila's text \cite{MP1}. Motivated by the above
papers, the authors in \cite{NB, NB1} introduced and studied a
multifractal formalism based on the Hewitt-Stromberg measures.
However, we point out that this formalism is completely parallel to
Olsen's multifractal formalism introduced in \cite{OL, OL1} which is
based on the Hausdorff and packing measures.

In the present paper we pursue those kinds of studies and consider
the  multifractal formalism developed in \cite{O}. We will start by
introducing the multifractal Hewitt-Stromberg measures and
dimensions which slightly differ from those introduced in \cite{NB,
NB1}. Our approach is to consider the behavior of multifractal
Hewitt-Stromberg functions dimensions under projections. These
dimension functions provide a natural mathematical framework within
which one can study multifractals since, for a wide class of
measures, they are related to multifractal spectra via the Legendre
transform. As an application, we study the multifractal analysis
based on Hewitt-Stromberg measures of the projections of a measure
and obtain general results for the multifractal analysis of the
orthogonal projections on $m$-dimensional linear subspaces of a
measure $\mu$ satisfying the multifractal formalism.

We will now give a brief description of the organization of the
paper. In the next section, we recall the definitions  of the
various multifractal dimensions and measures investigated in the
paper. Section \ref{sec1} recalls the multifractal formalism
introduced in \cite{O}. In Section \ref{sec2} we introduce the
multifractal Hewitt-Stromberg measures and separator functions
which slightly differ from those introduced in \cite{NB, NB1}, and
study their properties. Section \ref{Sec3} contain our main results.
The proofs are given in Sections \ref{Sec4}-\ref{Sec5}. The paper is
concluded with Section \ref{Sec6} which contains an application of
our main results.

\section{Preliminaries}

\subsection{Multifractal Hausdorff measure and packing
measure and dimensions}\label{sec1}

We start by recalling the multifractal formalism introduced by
O'Neil in \cite{O}. The key ideas behind the fine multifractal
formalism in \cite{O} are certain measures of Hausdorff-packing type
which are tailored to see only the multifractal decomposition sets
$E_\mu(\alpha)$. These measures are natural multifractal
generalizations of the centered Hausdorff measure and the packing
measure and are motivated by the $\tau_\mu$-function which appears
in the  multifractal formalism. We first recall the definition of
the multifractal Hausdorff measure and the the multifractal packing
measure. Let $\mu$ be a compactly supported probability measure on
$\mathbb{R}^n$. For $q, t\in\mathbb{R}$, $E \subseteq\supp\mu$ and
$\delta>0$, we define the multifractal packing pre-measure,

$$\overline{{\mathscr P}}^{q,t}_{\mu,\delta}(E) =\displaystyle
\sup \left\{\sum_i \mu\left(B\left(x_i,\frac{r_i}3\right)\right)^q
r_i^t\right\},$$ where the supremum is taken over all
$\delta$-packing of $E$,
$$\overline{{\mathscr P}}^{q,t}_{\mu}(E)
=\displaystyle\inf_{\delta>0}\overline{{\mathscr
P}}^{q,t}_{\mu,\delta}(E).$$

 \bigskip
The function $\overline{{\mathscr P}}^{q,t}_{\mu}$ is increasing but
not $\sigma$-subadditive. That is the reason why O'Neil introduced
the modification of the multifractal packing measure ${\mathscr
P}^{q,t}_{\mu}$:
$$
{\mathscr P}^{q,t}_{\mu}(E) = \inf_{E \subseteq \bigcup_{i}E_i}
\sum_i \overline{\mathscr P}^{q,t}_{\mu}(E_i).
 $$
In a similar way, we define the Hausdorff measure,
$${{\mathscr
H}}^{{q},t}_{{\mu},\delta}(E) = \displaystyle\inf \left\{\sum_i
\mu\Big(B(x_i,3r_i)\Big)^q ~r_i^t\;\Big|\;\Big(B(x_i,
r_i)\Big)_i\text{is a $\delta$-covering of}\; E\right\},
 $$
and
 $$
{{\mathscr H}}^{{q},t}_{{\mu}}(E) =
\displaystyle\sup_{\delta>0}{{\mathscr H}}^{q,t}_{{\mu},\delta}(E).
 $$

The functions ${\mathscr P}^{q,t}_{\mu}$ and ${\mathscr
H}^{q,t}_{\mu}$ are metric outer measures and thus measures on the
family of Borel subsets of $\mathbb{R}^n$. An important feature of
the pre-packing, packing and Hausdorff measure is that ${\mathscr
P}^{q,t}_{\mu}\leq{\overline{\mathscr P}}^{q,t}_{\mu}$ and there is
a constant $c$ depending also on the dimension of the ambient space,
such that ${\mathscr H}^{q,t}_{\mu}\leq c~{\mathscr P}^{q,t}_{\mu}$
(see \cite{O}).

\smallskip\smallskip
The functions $\overline{{\mathscr P}}^{q,t}_{\mu}$,  ${\mathscr
P}^{q,t}_{\mu}$ and ${\mathscr H}^{q,t}_{\mu}$ assign, in the usual
way, a dimension to each subset $E$ of $\supp\mu$. They are
respectively denoted by $\Lambda_{\mu}^q(E)$, $B_{\mu}^q(E)$ and
$b_{\mu}^q(E)$.

\par\noindent\begin{enumerate}\item  There exists a unique number
$\Lambda_{\mu}^q(E)\in[-\infty,+\infty]$ such that
 $$
\overline{{\mathscr P}}^{q,t}_{\mu}(E)=\left\{\begin{matrix}
\infty &\text{if}& t< \Lambda_{\mu}^q(E),\\ \\
0 & \text{if}&  \Lambda_{\mu}^q(E) < t.
\end{matrix}\right.
 $$
\item  There exists a unique number
$B_{\mu}^q(E)\in[-\infty,+\infty]$ such that
 $$
{\mathscr P}^{q,t}_{\mu}(E)=\left\{\begin{matrix}
\infty &\text{if}& t <B_{\mu}^q(E),\\
 \\
 0 & \text{if}&  B_{\mu}^q(E) < t.\end{matrix}\right.
 $$
 \item  There exists a unique number
$b_{\mu}^q(E)\in[-\infty,+\infty]$ such that
 $$
{\mathscr H}^{q,t}_{\mu}(E)=\left\{\begin{matrix}
\infty &\text{if}& t <b_{\mu}^q(E),\\
 \\
 0 & \text{if}&  b_{\mu}^q(E) < t.\end{matrix}\right.
 $$
\end{enumerate}
The number $b_{\mu}^q(E)$ is an obvious multifractal analogue of the
Hausdorff dimension $\dim_H(E)$ of $E$ whereas $B_{\mu}^q(E)$ and
$\Lambda_{\mu}^q(E)$ are obvious multifractal analogues of the
packing dimension $\dim_P(E)$ and the pre-packing dimension
$\Delta(E)$ of $E$ respectively. In fact, it follows immediately
from the definitions that
$$
\dim_H(E)=b_{\mu}^0(E),\;\;\;\dim_P(E)=B_{\mu}^0(E)\quad\text{and}\quad\Delta(E)=\Lambda_{\mu}^0(E).
$$

\bigskip
We note that for all $q\in\mathbb{R}$
$$
b_{\mu}^q(\emptyset)=B_{\mu}^q(\emptyset)=\Lambda_{\mu}^q(\emptyset)=-\infty,
$$
and if $\mu(E)=0$, then
$$
b_{\mu}^q(E)=B_{\mu}^q(E)=\Lambda_{\mu}^q(E)=-\infty\quad\text{for}\quad
q>0.
$$
Next, we define the separator functions $\Lambda_{\mu}$, $B_{\mu}$
and $b_{\mu}$ : $\mathbb{R}\rightarrow [-\infty, +\infty]$ by,
\begin{center}
$\Lambda_{\mu}(q)=\Lambda_{\mu}^q(\supp\mu)$,
$B_{\mu}(q)=B_{\mu}^q(\supp\mu)$ and
$b_{\mu}(q)=b_{\mu}^q(\supp\mu).$
\end{center}
It is well known that the functions $\Lambda_{\mu}$, $B_{\mu}$ and
$b_{\mu}$ are decreasing. The functions $\Lambda_{\mu}$, $B_{\mu}$
convex and satisfying $b_{\mu}\leq B_{\mu}\leq \Lambda_{\mu}.$

\begin{proposition}\cite{O}
Let $\mu$ be compactly supported probability measure on
$\mathbb{R}^n$. Then, we have
\begin{enumerate}
\item For $q<1$, $0\leq b_{\mu}(q)\leq B_{\mu}(q)\leq
\Lambda_{\mu}(q).$
\item $b_{\mu}(1)=B_{\mu}(1)=
\Lambda_{\mu}(1)=0.$
\item For $q>1$, $b_{\mu}(q)\leq B_{\mu}(q)\leq
\Lambda_{\mu}(q)\leq 0.$
\end{enumerate}
\end{proposition}
\begin{remark}
The multifractal Hausdorff and packing measures introduced by O'Neil
are different from those developed by Olsen \cite{OL}, although,
when $\mu$ satisfies a doubling condition, the multifractal measures
are equivalent.
\end{remark}
The reader is referred to as Olsen's classical text \cite{OL} (see
also \cite{O, SB3}) for an excellent and systematic discussion of
the multifractal Hausdorff and packing measures and dimensions.

%%%%%%%%%%%%%%%%%%%%%%%%%%%%%%%%%%%%%%%%%%%%%%%%%%%%%%%%%%%%%%%%%%%%%%%%%%%%%%
\subsection{Multifractal Hewitt-Stromberg measures and separator
functions}\label{sec2}
%%%%%%%%%%%%%%%%%%%%%%%%%%%%%%%%%%%%%%%%%%%%%%%%%%%%%%%%%%%%%%%%%%%%%%%%%%%%%%%%%%%

 Our main reason for modifying the definitions of Attia et al. in \cite{NB, NB1}  is to allow us
to prove results for non necessary doubling measures. One main cause
and motivation is the fact that such characteristics is not in fact
preserved under projections. In the following, we will set up, for
$q, t \in \mathbb{R}$ and a compactly supported probability measure
$\mu$ on $\mathbb{R}^n$, the lower and upper multifractal
Hewitt-Stromberg measures ${\mathsf{H}}_\mu^{q,t}$ and
${\mathsf{P}}_\mu^{q,t}$. For $E\subseteq \supp \mu$, the
pre-measure of $E$ is defined by
$$
 {\mathsf C}_\mu^{q,t}(E)= \limsup_{r\to0} M_{\mu,r}^q(E)~ r^t,
$$
where
$$
M_{\mu,r}^q(E)=\sup \left\{\displaystyle \sum_i
\mu\left(B\left(x_i,\frac
r3\right)\right)^q\;\Big|\;\Big(B(x_i,r)\Big)_i\;\text{is a packing
of}\; E \right\}.
$$
Observe that  ${\mathsf C}_\mu^{q,t}$ is increasing and ${\mathsf
C}_\mu^{q,t}(\emptyset ) =0$. However it is not $\sigma$-additive.
For this, we introduce  the ${\mathsf{P}}_\mu^{q,t}$-measure defined
by
$$
{\mathsf{P}}_\mu^{q,t}(E)=\inf \left\{\displaystyle \sum_i{\mathsf
C}_\mu^{q,t}(E_i)\;\Big|\; E\subseteq\bigcup_i E_i\; \text{and
the}\; E_i'\text{s are bounded} \right\}.
 $$
In a similar way we define
$$
{\mathsf L}_\mu^{q,t}(E)= \liminf_{r\to0} M_{\mu,r}^q(E)~ r^t.
$$
Since ${\mathsf L}_\mu^{q,t}$  is  not countably subadditive, one
needs a standard modification to get an outer measure. Hence, we
modify the definition as follows
 $$
{\mathsf{H}}_\mu^{q,t}(E)=\inf \left\{\displaystyle \sum_i{\mathsf
L}_\mu^{q,t}(E_i)\;\Big|\; E\subseteq\bigcup_i E_i\; \text{and
the}\; E_i'\text{s are bounded} \right\}
 $$

 \bigskip \bigskip
The measure $\mathsf{ H}^{q,t}_{\mu}$ is of course a multifractal
generalization of the lower $t$-dimensional Hewitt-Stromberg measure
${\mathsf{H}}^t$, whereas $\mathsf{ P}^{q,t}_{\mu}$ is a
multifractal generalization of the upper $t$-dimensional
Hewitt-Stromberg measures ${\mathsf{P}}^t$. In fact, it is easily
seen that, for $t>0$,  one has
$$
\mathsf{ H}^{0,t}_{\mu}={\mathsf{H}}^t\quad\text{and}\quad\mathsf{
P}^{0,t}_{\mu}={\mathsf{P}}^t.
$$

The following result describes some of the basic properties of the
multifractal Hewitt-Stromberg measures including the fact that
${\mathsf H}_\mu^{q,t}$ and ${\mathsf P}_\mu^{q,t}$ are Borel metric
outer measures and summarises the basic inequalities satisfied by
the multifractal Hewitt-Stromberg measures, the multifractal
Hausdorff measure and the multifractal packing measure.
\begin{theorem}\label{HHPP}
Let $q, t \in \mathbb{R}$. Then for every set $E\subseteq
\mathbb{R}^n$ we have
\begin{enumerate}
\item the set functions $\mathsf{H}_\mu^{q,t}$ and $\mathsf{P}_\mu^{q,t}$ are
metric outer measures and thus they are measures on the Borel
algebra.
\item There exists an integer $\xi\in\mathbb{N}$, such that $$
{\mathscr H}^{q,t}_{\mu}(E)\leq \xi\mathsf{H}_\mu^{q,t}(E)$$and $$
\begin{aligned}
&\mathsf{L}_\mu^{q,t}(E)\leq\mathsf{C}_\mu^{q,t}(E)\leq
\overline{{\mathscr P}}^{q,t}_{\mu}(E)\\&\;\;\;\mathrm{\vee
{l}}\quad\quad\quad\;\mathrm{\vee l}\quad\quad\quad\;\mathrm{\vee l}
\\&\mathsf{H}_\mu^{q,t}(E)\leq \mathsf{P}_\mu^{q,t}(E)\leq {\mathscr
P}^{q,t}_{\mu}(E).\end{aligned}$$
\end{enumerate}
\end{theorem}
The proof of the first part is straightforward and mimics that in
\cite[Theorem 2.1]{NB}. The proof of second part is a
straightforward application of Besicovitch's covering theorem and we
omit it here (we can see also \cite[Theorem 2.1]{NB}).

\bigskip\bigskip
The measures $\mathsf{ H}^{q,t}_{\mu}$ and $\mathsf{ P}^{q,t}_{\mu}$
and the pre-measures  $\mathsf{L}^{q,t}_{\mu}$ and
${\mathsf{C}}_{\mu}^{q,t}$ assign in the usual way a multifractal
dimension to each subset $E$ of $\mathbb{R}^n$, they are
respectively denoted by $\mathsf{b}_{\mu}^q(E)$,
$\mathsf{B}_{\mu}^q(E)$, ${\mathscr{L}}_{\mu}^q(E)$ and
$\mathsf{\Delta}_{\mu}^q(E)$.
\begin{proposition} Let $q \in \mathbb{R}$ and $E\subseteq
\mathbb{R}^n$. Then
\par\noindent\begin{enumerate}
\item  there exists a unique number ${\mathsf  b}_{\mu}^q(E)\in[-\infty,+\infty]$ such that
 $$
\mathsf{H}^{q,t}_{\mu}(E)=\left\{\begin{matrix}  \infty &\text{if}& t < {\mathsf  b}_{\mu}^q(E),\\
 \\
 0 & \text{if}&  {\mathsf  b}_{\mu}^q(E) < t,\end{matrix}\right.
 $$
\item  there exists a unique number ${\mathsf  B}_{\mu}^q(E)\in[-\infty,+\infty]$ such that
 $$
\mathsf{P}^{q,t}_{\mu}(E)=\left\{\begin{matrix}  \infty &\text{if}& t < {\mathsf  B}_{\mu}^q(E),\\
 \\
 0 & \text{if}&  {\mathsf  B}_{\mu}^q(E) < t,\end{matrix}\right.
 $$

\item  there exists a unique number ${\mathsf  \Delta}_{\mu}^q(E)\in[-\infty,+\infty]$ such that
 $$
\mathsf{C}^{q,t}_{\mu}(E)=\left\{\begin{matrix}  \infty &\text{if}& t < {\mathsf  \Delta}_{\mu}^q(E),\\
 \\
 0 & \text{if}&  {\mathsf  \Delta }_{\mu}^q(E) < t,\end{matrix}\right.
 $$
 \item there exists a unique number
$\mathscr{L}_{\mu}^q(E)\in[-\infty,+\infty]$ such that
 $$
{\mathsf L}^{q,t}_{\mu}(E)=\left\{\begin{matrix}  \infty &\text{if}& t < {\mathscr{L}}_{\mu}^q(E),\\
 \\
 0 & \text{if}&  {\mathscr{L}}_{\mu}^q(E) < t.\end{matrix}\right.
 $$
\end{enumerate}
In addition, we have
$${\mathsf  b}_{\mu}^q(E) \le   {\mathsf  B }_{\mu}^q(E) \le  {\mathsf  \Delta}_{\mu}^q(E).$$
\end{proposition}
The number ${\mathsf  b}_{\mu}^q(E)$ is an obvious multifractal
analogue of the lower Hewitt-Stromberg dimension
$\underline{\dim}_{MB}(E)$ of $E$ whereas ${\mathsf  B}_{\mu}^q(E)$
is an obvious multifractal analogues of the upper Hewitt-Stromberg
dimension $\overline{\dim}_{MB}(E)$ of $E$. Observe that the number
${\mathscr{L}}_{\mu}^q(E)$ is an obvious multifractal analogue of
the lower box-dimension $\underline{\dim}_{B}(E)$ of $E$ whereas
${\mathsf  \Delta}_{\mu}^q(E)$ is an obvious multifractal analogues
of the upper box-dimension $\overline{\dim}_{B}(E)$ of $E$. In fact,
it follows immediately from the definitions that
$$
{\mathscr{L}}_{\mu}^0(E)=\underline{\dim}_{B}(E),\;\;{\mathsf
\Delta}_{\mu}^0(E)=\overline{\dim}_{B}(E)$$ and $${\mathsf
b}_{\mu}^0(E)=\underline{\dim}_{MB}(E),\;\;{\mathsf
B}_{\mu}^0(E)=\overline{\dim}_{MB}(E).
$$
See \cite{NB1, JuMaMiOlSt, Olll2} for precise definitions of these
dimensions.
\begin{remark}\label{remark1}
It follows from Theorem \ref{HHPP} that
$$\begin{aligned}
&{\mathscr{L}}_{\mu}^q(E) \le  {\mathsf \Delta}_{\mu}^q(E)\leq
\Lambda_{\mu}^q(E)\\&\;\; \mathrm{\vee l}\;\qquad\quad\mathrm{\vee
l}\;\qquad\quad\mathrm{\vee l}
\\b_{\mu}^q(E) \le& \;{\mathsf b}_{\mu}^q(E) \le {\mathsf
B}_{\mu}^q(E) \le B_{\mu}^q(E).\end{aligned}$$
\end{remark}
The definition of these dimension functions makes it clear that they
are counterparts of the $\tau_\mu$-function which appears in the
{\it multifractal formalism}. This being the case, it is important
that they have the properties described by the physicists. The next
theorem shows that these functions do indeed have some of these
properties.
\begin{theorem}\label{th1}
Let $q\in\mathbb{R}$ and $E\subseteq \mathbb{R}^n$.
\begin{enumerate}
\item The functions $q \mapsto \mathsf{ H}^{q,t}_{\mu}(E)$,  $\mathsf{P}^{q,t}_{\mu}(E)$,  $ \mathsf{C}^{q,t}_{\mu}(E)$ are decreasing.
\item The functions $t \mapsto \mathsf{ H}^{q,t}_{\mu}(E)$,  $\mathsf{P}^{q,t}_{\mu}(E)$,  $ \mathsf{C}^{q,t}_{\mu}(E)$ are
decreasing.
\item The functions $q \mapsto { {\mathsf  b}}^{q}_{\mu}(E)$, $ {{\mathsf  B}}^{q}_{\mu}(E)$,  $ {{\mathsf  \Delta}}^{q}_{\mu}(E)$ are decreasing.
\item The functions $q \mapsto  {{\mathsf  B}}^{q}_{\mu}(E)$,  $ {{\mathsf  \Delta}}^{q}_{\mu}(E)$ are convex.
\end{enumerate}
\end{theorem}
The proof of this is straightforward and mimics that in
\cite[Theorem 3]{NB1}.

\bigskip\bigskip
We note that for all $q\in\mathbb{R}$
$$
{\mathsf b}_{\mu}^q(\emptyset)={\mathsf
B}_{\mu}^q(\emptyset)={\mathsf  \Delta}_{\mu}^q(\emptyset)=-\infty,
$$
and if $\mu(E)=0$, then
$$
{\mathsf
 b}_{\mu}^q(E)={\mathsf  B}_{\mu}^q(E)={\mathsf  \Delta}_{\mu}^q(E)=-\infty\quad\text{for}\quad
q>0.
$$
Next, we define the separator functions ${\mathsf  \Delta}_{\mu}$,
${\mathsf  B}_{\mu}$ and ${\mathsf  b}_{\mu}$ :
$\mathbb{R}\rightarrow [-\infty, +\infty]$ by,
\begin{center}
${\mathsf  \Delta}_{\mu}(q)={\mathsf  \Delta}_{\mu}^q(\supp\mu)$,
${\mathsf  B}_{\mu}(q)={\mathsf  B}_{\mu}^q(\supp\mu)$ and ${\mathsf
b}_{\mu}(q)={\mathsf  b}_{\mu}^q(\supp\mu).$
\end{center}
The multifractal formalism based on the measures ${\mathsf
H}^{q,t}_{\mu}$ and ${\mathsf P}^{q,t}_{\mu}$ and the dimension
functions $\mathsf{b}_{\mu}$, $\mathsf{B}_{\mu}$ and
$\mathsf{\Delta}_{\mu}$ provides a natural, unifying and very
general multifractal theory which includes all the hitherto
introduced multifractal parameters, i.e., the multifractal spectra
functions $\alpha\mapsto
\mathsf{f}_\mu(\alpha)=:\underline{\dim}_{MB}E_\mu(\alpha)$ and
$\alpha\mapsto
\mathsf{F}_\mu(\alpha)=:\overline{\dim}_{MB}E_\mu(\alpha)$, the
multifractal box dimensions. The dimension functions
$\mathsf{b}_{\mu}$ and $\mathsf{B}_{\mu}$ are intimately related to
the spectra functions $\mathsf{f}_\mu$ and $\mathsf{F}_\mu$ (see
\cite{NB1}), whereas the dimension function $\mathsf{\Delta}_{\mu}$
is closely related to the upper box spectrum (more precisely, to the
upper multifractal box dimension function $\overline{\tau}_{\mu}$,
see Proposition \ref{newbox}).

%%%%%%%%%%%%%%%%%%%%%%%%%%%%%%%%%%%%%%%%%%%%%%%%%%%%%%%%%%%%%%%%%%%%%%%%%%%%%%%ù
\section{Main results}\label{Sec3}

Let $\mu$ be a compactly supported probability measure on
$\mathbb{R}^n$ and $q\in \mathbb{R}$. In the following, we require
an alternative characterization of the upper and lower multifractal
box-counting dimensions of $\mu$  in terms of a potential obtained
by convolving $\mu$ with a certain kernel. For this purpose let us
introduce some notations. For $1\leq s\leq n$ and $r>0$ we define
 $$
\begin{array}{llll}\label{j}
\phi_r^s: & \mathbb{R}^n & \longrightarrow & \mathbb{R} \\
& x & \longmapsto & \min\Big\{1,\: r^s|x|^{-s}\Big\},
\end{array}
$$
and
\begin{equation*}
\label{r}\mu\ast\phi_r^s(x)=\int\min\Big\{1,\:
r^s|x-y|^{-s}\Big\}d\mu(y).\end{equation*} Let $E$ be a compact
subset of $\supp\mu$. For $1\leq s\leq n$ and $q>1$, write
$$
N_{\mu,r}^{q,s}(E)=
\int_E\Big(\mu\ast\phi_{r/3}^s(x)\Big)^{q-1}d\mu(x),
 $$
and
$$
\overline{\tau}_{\mu}^{q,s}(E)=\displaystyle\limsup_{r\to0}
\frac{\log N_{\mu,r}^{q,s}(E)}{-\log r}\quad\text{and}\quad
\underline{\tau}_{\mu}^{q,s}(E)=\displaystyle\liminf_{r\to0}
\frac{\log N_{\mu,r}^{q,s}(E)}{-\log r}.
$$

These definitions are, frankly, messy, indirect and unappealing. In
an attempt to make the concept more attractive, we present here an
alternative approach to the dimensions
$\underline{\tau}_{\mu}^{q,s}$ and  $\overline{\tau}_{\mu}^{q,s}$,
and their applications to projections in terms of a potential
obtained by convolving $\mu$ with a certain kernel. For $E$ a
compact subset of $\supp\mu$ we can try to decompose $E$ into a
countable number of pieces $E_1, E_2, . . .$ in such a way that the
largest piece has as small a dimension as possible. The present
approach was first used by Falconer in \cite[Section 3.3]{Falconer}
and further developed by O'Neil and Selmi in \cite{O, SB3}. This
idea leads to the following modified dimensions in terms of the
convolutions:
\begin{eqnarray*}
\underline{\mathfrak{T}}_\mu^{q,s}(E)&=&\inf\left\{\sup_{1\leq
i<\infty}\underline{\tau}_{\mu}^{q,s}(E_i)\;\Big|\;
E\subseteq\bigcup_i E_i\;\;\text{with each $E_i$ compact }\right\},
\end{eqnarray*}
\begin{eqnarray*}
\overline{\mathfrak{T}}_\mu^{q,s}(E)&=&\inf\left\{\sup_{1\leq
i<\infty}\overline{\tau}_{\mu}^{q,s}(E_i)\;\Big|\;
E\subseteq\bigcup_i E_i\;\;\text{with each $E_i$ compact }\right\}
\end{eqnarray*}
and
$$
\underline{\mathfrak{T}}_\mu^{s}(q)=\underline{\mathfrak{T}}_\mu^{q,s}(\supp\mu)
\quad\text{and}\quad\overline{\mathfrak{T}}_\mu^{s}(q)=\overline{\mathfrak{T}}_\mu^{q,s}(\supp\mu)\quad\text{for
all}\quad s\geq 1.
$$

\bigskip
Let $m$ be an integer with $0<m\leq n$ and $G_{n,m}$ the
Grassmannian manifold of all $m$-dimensional linear subspaces of
$\mathbb{R}^n$. Denote by $\gamma_{n,m}$ the invariant Haar measure
on $G_{n,m}$ such that $\gamma_{n,m}(G_{n,m})=1$. For $V\in
G_{n,m}$, we define the projection map $\pi_V:
\mathbb{R}^n\longrightarrow V$ as the usual orthogonal projection
onto $V$. Then, the set $\{\pi_V,\; V \in G_{n,m}\}$ is compact in
the space of all linear maps from $\mathbb{R}^n$ to $\mathbb{R}^m$
and the identification of $V$ with $\pi_V$ induces a compact
topology for $G_{n,m}$. Also, for a Borel probability measure $\mu$
with compact support $supp\mu \subset\mathbb{R}^n$ and for $V\in
G_{n,m}$, we denote by $\mu_V,$ the projection of $\mu$ onto $V$,
i.e.,
 $$
\mu_V(A)=\mu\circ\pi_V^{-1}(A)\quad \forall A\subseteq V.
 $$

Since $\mu$ is compactly supported and $supp\mu_V=\pi_V(supp\mu)$
for all $V\in G_{n,m}$, then, for any continuous function $f:
V\longrightarrow\mathbb {R}$, we have
 $$
\displaystyle\int_V fd\mu_V=\int f(\pi_V(x))d\mu(x),
 $$
whenever these integrals exist. Then for all $V\in G_{n,m}$, $x\in
\mathbb{R}^n$ and $0<r<1$, we have
\begin{equation*}
\label{r}\mu\ast\phi_r^m(x)=\int\mu_V(B(x_V,r))dV=\int\min\Big\{1,\:
r^m|x-y|^{-m}\Big\}d\mu(y).\end{equation*}

\bigskip
In the following we compare the lower and upper multifractal
Hewitt-Stromberg dimensions of a set $E$ of $\mathbb{R}^n$ with
respect to a measure $\mu$ with those of their projections onto
$m$-dimensional subspaces.
\begin{theorem}\label{TH1}
Let $\mu$ be a compactly supported probability measure on
$\mathbb{R}^n$ and $E\subseteq\supp\mu$. For $q\leq 1$ and all $V\in
G_{n,m}$, we have
$$
{\mathscr{L}}_{\mu_V}^{q}(\pi_V(E))\leq{\mathscr{L}}_{\mu}^{q}(E),\;\;\;
{\mathsf \Delta}_{\mu_V}^{q}(\pi_V(E))\leq {\mathsf
\Delta}_{\mu}^{q}(E)
$$
and
$$
{\mathsf b}_{\mu_V}^{q}(\pi_V(E))\leq {\mathsf
b}_{\mu}^{q}(E),\;\;\; {\mathsf B}_{\mu_V}^{q}(\pi_V(E))\leq
{\mathsf B}_{\mu}^{q}(E).
$$
\end{theorem}

\bigskip
The following result presents alternative expressions of the
multifractal dimension functions ${\mathscr{L}}_{\mu}^{q}(E)$ and
${\mathsf \Delta}_{\mu}^{q}(E)$ of a set $E$ and that of its
orthogonal projections.
\begin{theorem} \label{THLD} Let $E$ be a compact subset of $\supp\mu$. Then, we have
\begin{enumerate}
\item for all $q>1$ and $V\in G_{n,m},$
 $$
{\mathscr{L}}_{\mu_V}^{q}(\pi_V(E))\geq\underline{\tau}_{\mu}^{q,m}(E)\geq
{\mathscr{L}}_{\mu}^{q}(E)
$$
and
$$
{\Delta}_{\mu_V}^{q}(\pi_V(E))\geq\overline{\tau}_{\mu}^{q,m}(E)\geq
{\mathsf \Delta}_{\mu}^{q}(E).
 $$

\item For all $1<q\leq2$ and $\gamma_{n,m}$-almost every  $V\in
G_{n,m},$
 $$
{\mathsf\Delta}_{\mu_V}^{q}(\pi_V(E))=\overline{\tau}_{\mu}^{q,m}(E)=\max\Big(m(1-q),
{\mathsf  \Delta}_{\mu}^{q}(E)\Big)
 $$
and
 $$
{\mathscr{L}}_{\mu_V}^{q}(\pi_V(E))=\underline{\tau}_{\mu}^{q,m}(E).
 $$
 \item For all $q>2$ and $\gamma_{n,m}$-almost every  $V\in
G_{n,m},$
\begin{enumerate}
 \bigskip
\item If ${\mathsf  \Delta}_{\mu}^{q}(E)\geq-m$ then
${\mathsf
\Delta}_{\mu_V}^{q}(\pi_V(E))=\overline{\tau}_{\mu}^{q,m}(E)=
{\mathsf  \Delta}_{\mu}^{q}(E).$
 \bigskip
\item $
{\mathscr{L}}_{\mu_V}^{q}(\pi_V(E))=\max\Big(
m(1-q),\underline{\tau}_{\mu}^{q,m}(E)\Big).$
\end{enumerate}
\end{enumerate}
\end{theorem}

\bigskip
In Theorem \ref{TH2} we show that ${\mathsf B}_\mu^q(E)$ is
preserved under $\gamma_{n,m}$-almost every orthogonal projection
for $q>1$.
\begin{theorem} \label{TH2} Let $E$ be a compact subset of $\supp\mu$ and $q>1$.
\begin{enumerate}
\item For all $V\in G_{n,m}$, we have
$$
{\mathsf B}_{\mu_V}^{q}(\pi_V(E))\geq {\mathsf B}_{\mu}^{q}(E).
$$
\item If $1<q\leq2$, one has
$$
{\mathsf
B}_{\mu_V}^{q}(\pi_V(E))=\overline{\mathfrak{T}}_{\mu}^{q,m}(E)=\max\Big(m(1-q),
{\mathsf B}_{\mu}^{q}(E)\Big),\;\text{for $\gamma_{n,m}$-almost
every}\;\; V\in G_{n,m}.
$$
 \item If $q>2$ and $({E}_i)_i$  is a cover
of $E$ by a countable collection of compact sets is such that
${\mathsf  \Delta}_{\mu}^q(E_i)\geq -m$ for all $i$, then
 $$
{\mathsf
B}_{\mu_V}^{q}(\pi_V(E))=\overline{\mathfrak{T}}_{\mu}^{q,m}(E)={\mathsf
B}_{\mu}^{q}(E),\;\text{for $\gamma_{n,m}$-almost every}\;\; V\in
G_{n,m}.
$$
\end{enumerate}
\end{theorem}

\bigskip
The next theorem enables us to study the lower multifractal
Hewitt-Stromberg dimension of the projection of sets on
$m$-dimensional linear subspaces for $q>1$. In particular, we prove
that ${\mathsf b}_\mu^q(E)$ is not preserved under $\gamma
_{n,m}$-almost every orthogonal projection for $q>1$.
\begin{theorem} \label{TH3} Let $E$ be a compact subset of $\supp\mu$ and $q>1$.
\begin{enumerate}
\item For all $V\in G_{n,m}$, we have
$$
{\mathsf b}_{\mu_V}^{q}(\pi_V(E))\geq {\mathsf b}_{\mu}^{q}(E).
$$
\item If $1<q\leq2$, one has
$$
{\mathsf
b}_{\mu_V}^{q}(\pi_V(E))=\underline{\mathfrak{T}}_{\mu}^{q,m}(E),\;\text{for
$\gamma_{n,m}$-almost every}\;\; V\in G_{n,m}.
$$
 \item If $q>2$, then
 $$
{\mathsf b}_{\mu_V}^{q}(\pi_V(E))=\max\Big(m(1-q),
\underline{\mathfrak{T}}_{\mu}^{q,m}(E)\Big),\;\text{for
$\gamma_{n,m}$-almost every}\;\; V\in G_{n,m}.
$$
\end{enumerate}
\end{theorem}

\bigskip\bigskip
For an integer $s$ with $1\leq m \leq s<n$, we define the $s$-energy
of a measure $\mu$ by
\begin{eqnarray*}\label{energ}
I_s(\mu)=\displaystyle\int\int|x-y|^{-s}d\mu(x)d\mu(y).
\end{eqnarray*}
Frostman \cite{FR} showed that the Hausdorff dimension of a Borel
subset $E$ of $\mathbb{R}^n$ is the supremum of the positive reals
$s$ for which there exists a Borel probability measure $\mu$
charging $E$ and for which the $s$-energy of $\mu$ is finite. This
characterization is used by Kaufmann \cite{K} and Mattila \cite{MP1}
to prove their results on the preservation of the Hausdorff
dimension. The condition $I_s(\mu)<\infty$ implies that
$\dim_H(\mu)\geq s$. On the other hand, if $\mu(B(x,r))\leq r^s$,
for all $x$ and all sufficiently small $r$ then $\mu$ has a finite
$s$-energy. Notice that Mattila \cite{MP1} proved that if $I_m(\mu)$
is finite, then for almost every $m$-dimensional subspace $V$, the
measure $\mu_V$ is absolutely continuous with respect to Lebesgue
measure $\mathcal{L}_V^m$ on $V$ identified with $\mathbb{R}^m$ and
$\mu_V\in L^2(V)$, where $\mathcal{L}_V^m(E)=\mathcal{L}^m(E\cap V)$
for $E\subset\mathbb{R}^m$. The following result enables us to
describe the behavior of large measures under projection.
\begin{theorem}\label{THbB}
 Suppose that $\mu$  is a compactly supported Radon measure on  $\mathbb{R}^{n}$  and $0<m \leq s<n$
 are such  that  $I_{s}(\mu)<\infty.$ Then
\begin{enumerate}
\item if  $2 m<s<n$,  then $\text{for
$\gamma_{n,m}$-almost every}\;\; V\in G_{n,m}$  and  $q \geq 0$
$$ {\mathsf b}_{\mu_{V}}({q})={\mathsf B}_{\mu_V}({q})=m(1-q),$$

\item if  $m \leq s \leq 2 m$, then $\text{for
$\gamma_{n,m}$-almost every}\;\; V\in G_{n,m}$  and $q \geq0$
$$
 m(1-q) \leq {\mathsf b}_{\mu_V}(q)\leq {\mathsf B}_{\mu_V}(q) \leq \max \Big(m(1-q),-\frac{s q} 2\Big).$$
\end{enumerate}
\end{theorem}
\begin{remark}
Fix $0<m\leq n$ and let $\mu$ be a self-similar measure on
$\mathbb{R}^n$ with support equal to $K$ such that $\dim_P(K)=s\leq
m$. Let $q\geq0$ and $({E}_i)_i$  be a cover of $E$ by a countable
collection of compact sets is such that ${\mathsf
\Delta}_{\mu}^q(E_i)\geq -m$ for all $i$. By using Theorems
\ref{TH2} and \ref{TH3} and   \cite[Corollary 5.12]{O}, we have for
$\gamma_{n,m}$-almost every $V\in G_{n,m}$
\begin{eqnarray*}
&B_{{\mu}_V}({q})={\mathsf b}_{{\mu}_V}({q})={\mathsf
B}_{{\mu}_V}({q})=b_{{\mu}_V}({q})= b_{{\mu}}({{q}})=\\&{\mathsf
b}_{{\mu}}({{q}})={\mathsf
B}_{{\mu}}({{q}})=B_{{\mu}}({q})=\underline{\mathfrak{T}}_{\mu}^{m}(q)=\overline{\mathfrak{T}}_{\mu}^{m}(q).
\end{eqnarray*}
\end{remark}

\section{Proof of Theorem \ref{TH1}}\label{Sec4} When $q<0$ it suffice to observe that if $\Big(B\left(x_{i}, r\right)\Big)_{i
\in \mathbb{N}}$  is a centered packing of  $\pi_{V}(E)$  then
$\Big(B\left(y_{i}, r\right)\Big)_{i \in \mathbb{N}}$  is a centered
packing of  $E$ where $y_{i} \in E$ is such that
$x_{i}=\pi_{V}\left(y_{i}\right)$ which implies that  $\mu_{V}
(B\left(x_{i} \frac r3\right))^{q} \leq \mu (B\left(y_{i}, \frac
r3\right))^{q}$.  This easily gives the desired result.

\bigskip\bigskip
\noindent For $0\leq q\leq1,$ fix $V\in G_{n,m}$ and let
$\Big(B(x_i,r)\Big)_{i \in \mathbb{N}}$ be a packing of $\pi_V(E)$.
For each $i$ we consider the collection of balls as follows
 $$
\left\{B\left(y,\frac{r}9\right)\;\Big|\; y \in E \cap
\pi_V^{-1}\left (V\cap B\left(x_i,\frac{r}3\right)\right)\right\},
 $$
and let $E_i=E \cap \pi_V^{-1}\big (V\cap B(x_i,\frac{r}3)\big)$. So
Besicovitch's covering theorem (see \cite[Theorem 2.2]{O}) provides
a positive integer $\xi=\xi(n)$ and index sets, $
I_{i1},I_{i2},...,I_{i\xi}$ such that $ E_i\subset
\bigcup_{j=1}^{\xi}\bigcup_{i\in I_j }B\left(y,\frac{r}9\right)$ and
the subset $ \Big\{B(y,\frac{r}9)\;\big|\; y\in I_{ij}\Big\}$ is a
disjoint family for each $j$. Also, we observe that for a fixed $i$
and $j$, we may make a simple volume estimate to further subdivide
$J_{ij}$ into at most $7^n$ disjoint subfamilies, $J_{ij1}, J_{ij2},
. . . , J_{ij7^n}$ such that for each $j$ and $k$, $
\left\{B\left(y,\frac{r}3\right)\;\big|\; y\in \cup_i
I_{ijk}\right\}$ is a centered packing of $E$. Since $0<q\leq1$, we
have
\begin{eqnarray*}
\sum_i\mu_V\left(B\left(x_i,\frac {r} 3\right)\right)^q&\leq&
\sum_i\mu \left
(\bigcup_{j=1}^{\xi}\bigcup_{k=1}^{7^n}\bigcup_{y\in I_{ijk}} B\left(y,\frac {r} 9\right) \right )^q\\
&\leq& \sum_i\sum_{j=1}^{\xi}\sum_{k=1}^{7^n}\sum_{y\in I_{ijk}}\mu
\left ( B\left(y,\frac {r} 9\right)\right)^q\\
&\leq& \xi 7^n M_{\mu,\frac r3}^{q}(E).
\end{eqnarray*}
This implies that
 $$
M_{\mu_V,r}^{q}(\pi_V(E))\leq
 \xi 7^n M_{\mu,\frac r
3}^{q}(E).
 $$
Letting $r\downarrow0$, now yields
\begin{equation*}
{\mathsf L}_{\mu_V}^{q,t}(\pi_V(E))\leq 3^t \xi 7^n {\mathsf
L}_{\mu}^{q,t}(E)\quad\text{and}\quad {\mathsf
C}_{\mu_V}^{q,t}(\pi_V(E))\leq 3^t \xi 7^n {\mathsf
C}_{\mu}^{q,t}(E).
\end{equation*}
We deduce from the previous inequalities that
$$
{\mathscr{L}}_{\mu_V}^{q}(\pi_V(E))\leq{\mathscr{L}}_{\mu}^{q}(E)\quad\text{and}\quad
{\mathsf \Delta}_{\mu_V}^{q}(\pi_V(E))\leq {\mathsf
\Delta}_{\mu}^{q}(E).
$$

\bigskip
Now, let $t>{\mathsf B}_{\mu}^{q}(E)$ which implies that
$\mathsf{P}_{\mu}^{q,t}(E)<\infty$, then we can choose $(E_i)_i$ a
covering of $E$ such that $\displaystyle\sum_i{\mathsf
C}_{\mu}^{q,t}(E_i)<1.$ We therefore conclude that
$\pi_V(E)\subseteq\displaystyle\bigcup_i\pi_V(E_i)$ and
\begin{eqnarray*}
\mathsf{P}_{\mu_V}^{q,t}(\pi_V(E)) &\leq& \displaystyle\sum_i
{\mathsf C}_{\mu_V}^{q,t}(\pi_V(E_i)) \leq 3^t\xi 7^n \sum_i{\mathsf
C}_{\mu}^{q,t}(E_i) \leq3^t\xi 7^n<\infty.
\end{eqnarray*}
This implies that
$$
{\mathsf B}_{\mu_V}^{q}(\pi_V(E))\leq t\quad\text{for all}\quad
t>{\mathsf B}_{\mu}^{q}(E).
$$
We now infer that
\begin{eqnarray}\label{ProjectionbB}
{\mathsf B}_{\mu_V}^{q}(\pi_V(E))\leq {\mathsf B}_{\mu}^{q}(E).
\end{eqnarray}
The proof of the statement ${\mathsf b}_{\mu_V}^{q}(\pi_V(E))\leq
{\mathsf b}_{\mu}^{q}(E)$ is identical to the proof of the statement
\eqref{ProjectionbB} and is therefore omitted.
\section{Proof of Theorems \ref{THLD}, \ref{TH2}, \ref{TH3} and \ref{THbB}}\label{Sec5}
We present the tools, as well as the intermediate results, which
will be used in the proof of our main results.
\subsection{Preliminary results}
Let $\mu$ be a compactly supported probability measure on
$\mathbb{R}^n$ and $q\in \mathbb{R}$. Recall that the upper and
lower multifractal box-counting dimensions $\overline{\tau}_{\mu}^q$
and $\underline{\tau}_{\mu}^q$ of $E$ are defined respectively  by
$$
\overline{\tau}_{\mu}^q(E)=\displaystyle\limsup_{r\to0} \frac{\log
M_{\mu,r}^q(E)}{-\log
r}\quad\text{and}\quad\underline{\tau}_{\mu}^q(E)=\displaystyle\liminf_{r\to0}
\frac{\log M_{\mu,r}^q(E)}{-\log r}.
$$
For technical convenience we shall assume that,, if $r\in (0,1)$:
$\frac{\log 0}{-\log r}=-\infty$.

\smallskip\smallskip

The next result is essentially a restatement of \cite[Proposition
4.2]{BB} and \cite[Proposition 5.1]{B} (see also \cite[Lemma 2.6
(a)]{FO} and \cite{SS}), and has recently been obtained in
\cite[Proposition 4.2]{SB3}.
\begin{proposition}\label{P1}
Let $E$ be a compact subset of $\supp\mu$. For $q>1$, we have
$$
\underline{\tau}_{\mu}^{q}(E)=\liminf_{r\to0}\frac{1}{-\log
r}\log\int_E\mu\left(B\left(x,\frac r3\right)\right)^{q-1}d\mu(x)
$$
and
$$
\overline{\tau}_{\mu}^{q}(E)=\limsup_{r\to0}\frac{1}{-\log
r}\log\int_E\mu\left(B\left(x,\frac r3\right)\right)^{q-1}d\mu(x).
$$
\end{proposition}

\bigskip\bigskip
 In the next we  investigate the relation between the lower and
upper multifractal Hewitt-Stromberg functions $\mathsf{b}_\mu$ and
$\mathsf{B}_\mu$ and the multifractal box dimension, the
multifractal packing dimension and the multifractal pre-packing
dimension.

\begin{proposition}\label{newbox}
Let $q\in\mathbb{R}$ and $\mu$ be a compact supported Borel
probability measure on $\mathbb{R}^n$. Then for every $E \subseteq
\supp\mu$ we have
$$
{\mathscr{L}}_{\mu}^q(E) = \underline{\tau}_{\mu}^q(E)\qquad
\text{and }\qquad {\mathsf \Delta}_{\mu}^q(E) =
\overline{\tau}_{\mu}^q(E)=\Lambda_{\mu}^q(E).
$$
\end{proposition}
\noindent{\bf Proof.} We will prove the first equality, the second
one is similar. Suppose that $$ \underline{\tau}_{\mu}^q(E) >
{\mathscr{L}}_{\mu}^q(E) +\epsilon\quad\text{ for some}\quad\epsilon
>0.$$ Then we can find $\delta >0$ such that for any $r\le \delta$,
$$
M_{\mu, r}^q (E)~ r^{{  \mathscr{L}}_{\mu}^q(E)  +\epsilon}
>1\quad\text{and then} \quad \mathsf{L}_{\mu}^{q, {\mathscr{L}}_{\mu}^q(E)
+\epsilon} \ge 1
$$
which is a contradiction. We therefore infer $$
\underline{\tau}_{\mu}^q(E) \le {\mathscr{L}}_{\mu}^q(E)
+\epsilon\;\;\text{ for any}\;\; \epsilon
>0.$$ The proof of the following statement  $$
\underline{\tau}_{\mu}^q(E) \ge {\mathscr{L}}_{\mu}^q(E) -
\epsilon\;\;\text{ for any}\;\;\epsilon
>0$$
is identical to the proof of the above statement and is therefore
omitted.

\bigskip
We have the following additional property.
\begin{proposition}\label{modbox}
Let $q\in\mathbb{R}$ and $\mu$ be a compact supported Borel
probability measure on $\mathbb{R}^n$. Then for every $E \subseteq
\supp \mu$ we have
$$
{\mathsf  b}_{\mu}^q(E) =\inf \left\{ \sup_{i}
{\mathscr{L}}_{\mu}^q(E_i) \;\Big|\;\displaystyle E\subseteq
\bigcup_i E_i, \;\; E_i \;\; \text{are bounded in } \;\;
\mathbb{R}^n\right\}
$$
and
$$
{\mathsf  B}_{\mu}^q(E) = \inf \left\{ \sup_{i}  {\mathsf
\Delta}_{\mu}^q (E_i) \;\Big|\;\displaystyle E\subseteq \bigcup_i
E_i, \;\; E_i \;\; \text{are bounded in } \;\; \mathbb{R}^n\right\}.
$$

\end{proposition}
\noindent{\bf Proof.} Denote
$$\beta = \inf \left\{ \sup_{i} {\mathscr{L}}_{\mu}^q(E_i)
\;\Big|\;\displaystyle E\subseteq \bigcup_i E_i, \;\; E_i \;\;
\text{are bounded in } \;\; \mathbb{R}^n\right\}.$$ Assume that
$\beta < {\mathsf  b}_{\mu}^q(E)$ and take $\alpha \in (\beta,
{\mathsf b}_{\mu}^q(E) )$. Then we can choose $\{E_i\}$ of bounded
subset of $E$ such that  $E\subseteq \cup_i E_i$, and $\sup_i
{\mathscr{L}}_{\mu}^q(E_i) < \alpha$. Observe that
$\mathsf{L}_{\mu}^{q, \alpha}(E_i) = 0$ which implies that
${{\mathsf H}}_{\mu}^{q, \alpha}(E) = 0$. It is a contradiction. Now
suppose that $ {\mathsf b}_{\mu}^q(E) < \beta$, then, for any
$\alpha \in ( {\mathsf b}_{\mu}^q(E), \beta ),$ we have ${\mathsf
H}_{\mu}^{q, \alpha}(E) = 0$. Thus, there exists $\{E_i\}$ of
bounded subset of $E$ such that $E\subseteq \cup_i E_i$, and $\sup_i
{\mathsf L}_{\mu}^{q, \alpha} (E_i) < \infty$. We conclude that,
$\sup_i {\mathscr{L}}_{\mu}^{q} (E_i) \le \alpha$. It is also a
contradiction. The proof of the second statement is identical to the
proof of the statement in the first part and is therefore omitted.

 \bigskip
\noindent The following proposition is a consequence of Propositions
\ref{newbox} and \ref{modbox}.
\begin{proposition}\label{ah} Let $E$ be a compact subset of $\supp\mu$ and $q\in
\mathbb{R}$. One has
\begin{eqnarray*}
\mathsf{b}_\mu^q(E)=\inf\left\{\sup_{1\leq
i<\infty}\underline{\tau}_{\mu}^{q}(E_i)\;\Big|\;
E\subseteq\bigcup_i E_i\;\;\text{with each $E_i$ compact}\;\right\}
\end{eqnarray*}
and
\begin{eqnarray*}
\mathsf{B}_\mu^q(E)&=&\inf\left\{\sup_{1\leq
i<\infty}\overline{\tau}_{\mu}^{q}(E_i)\;\Big|\; E\subseteq\bigcup_i
E_i\;\;\text{with each $E_i$ compact}\;\right\}\\
&=&\inf\left\{\sup_{1\leq i<\infty}\Lambda_{\mu}^{q}(E_i)\;\Big|\;
E\subseteq\bigcup_i E_i\;\;\text{with each $E_i$ compact}\;\right\}.
\end{eqnarray*}
\end{proposition}

\begin{proposition}\label{L1} Let $E$ be a subset of $\supp\mu$ and $q\in
\mathbb{R}$. Then we have
$$
B_\mu^q(E)= \mathsf{B}_\mu^q(E).
$$
\end{proposition}
\noindent{\bf Proof.} It follows immediately from Proposition
\ref{ah} and \cite[Proposition 4.1]{SB3}.

\bigskip
The following result presents alternative expressions of the upper
and lower multifractal box-counting dimensions in terms of the
convolutions as well as general relations between the  upper and
lower multifractal box-counting dimensions of a measure and that of
its orthogonal projections. This result has recently been obtained
in \cite[Theorem 4.1]{SB3}.
\begin{theorem} \label{TH222} Let $E$ be a compact subset of $\supp\mu$. Then, we have
\begin{enumerate}
\item for all $q>1$ and $V\in G_{n,m},$
 $$
\underline{\tau}_{\mu_V}^{q}(\pi_V(E))\geq\underline{\tau}_{\mu}^{q,m}(E)\quad\text{and}
\quad\overline{\tau}_{\mu_V}^{q}(\pi_V(E))\geq\overline{\tau}_{\mu}^{q,m}(E).
 $$

\item For all $1<q\leq2$ and $\gamma_{n,m}$-almost every  $V\in
G_{n,m},$
 $$
\overline{\tau}_{\mu_V}^{q}(\pi_V(E))=\overline{\tau}_{\mu}^{q,m}(E)=\max\Big(m(1-q),
\overline{\tau}_{\mu}^{q}(E)\Big)
 $$
and
 $$
\underline{\tau}_{\mu_V}^{q}(\pi_V(E))=\underline{\tau}_{\mu}^{q,m}(E).
 $$
 \item For all $q>2$ and $\gamma_{n,m}$-almost every  $V\in
G_{n,m},$
\begin{enumerate}
 \bigskip
\item If $-m\leq\overline{\tau}_{\mu}^{q}(E)$ then
$\overline{\tau}_{\mu_V}^{q}(\pi_V(E))=\overline{\tau}_{\mu}^{q,m}(E)=
\overline{\tau}_{\mu}^{q}(E).$
 \bigskip
\item $
\underline{\tau}_{\mu_V}^{q}(\pi_V(E))=\max\Big(
m(1-q),\underline{\tau}_{\mu}^{q,m}(E)\Big).$
\end{enumerate}
\end{enumerate}
\end{theorem}
 The assertion (2) is essentially a restatement
of the main result of Hunt et al. in \cite{HK} and  Falconer et al.
in \cite[Theorem 3.9]{FO}. The assertion (3) extends the result of
Hunt and Kaloshin (of Falconer and O'Neil) to the case $q> 2$
untreated in their work.

\subsection{ Proof of Theorem \ref{THLD}}
 Follows directly from Proposition \ref{newbox} and Theorem \ref{TH222}.
\subsection{ Proof of Theorem \ref{TH2}}
The proof of the first part of Theorem \ref{TH2} follows immediately
from Theorem \ref{TH222} (1.), Proposition \ref{modbox} and since
$\overline{\tau}_{\mu}^{q,m}(E)\geq\overline{\tau}_{\mu}^{q}(E)$. By
using Proposition \ref{L1}, then the proof of the second and the
third part of Theorem \ref{TH2} has recently been obtained in
\cite[Theorem 3.2]{SB3}.

\subsection{ Proof of Theorem \ref{TH3}}
\begin{enumerate}
\item  The proof of the first part of Theorem \ref{TH3} follows immediately
from Theorem \ref{TH222} (1.), Proposition \ref{modbox} and since
$\underline{\tau}_{\mu}^{q,m}(E)\geq\underline{\tau}_{\mu}^{q}(E)$.
\item If $s>\underline{\mathfrak{T}}_{\mu}^{q,m}(E)$ we may cover $E$ by a countable
collection of sets $E_i$, which we may take to be compact, such that
$\underline{\tau}_{\mu}^{q,m}(E_i)<s$. By using Theorem \ref{TH222}
(2.), we have $\underline{\tau}_{\mu_V}^{q}(\pi_V(E_i))\leq s$ for
$\gamma_{n,m}$-almost every  $V\in G_{n,m}.$ Proposition \ref{ah}
implies that ${\mathsf b}_{\mu_V}^{q}(\pi_V(E))\leq s$ for
$\gamma_{n,m}$-almost every $V\in G_{n,m}$ and so, ${\mathsf
b}_{\mu_V}^{q}(\pi_V(E))\leq
\underline{\mathfrak{T}}_{\mu}^{q,m}(E)$ for $\gamma_{n,m}$-almost
every  $V\in G_{n,m}.$

Now, if $s<\underline{\mathfrak{T}}_{\mu}^{q,m}(E)$. Fix $V\in
G_{n,m}$ and let $(\widetilde{E}_i)_i$ be a cover of the compact set
$\pi_V(E)$ by a countable collection of compact sets. Put for each
$i$, $E_i=E\cap\pi_V^{-1}(\widetilde{E}_i),$ then
$\sup_i\underline{\tau}_{\mu}^{q,m}(E_i)>s$. By using Theorem
\ref{TH222} (1.), we have
$\sup_i\underline{\tau}_{\mu_V}^{q}(\pi_V(E_i))\geq s$ and
$\sup_i\underline{\tau}_{\mu_V}^{q}(\widetilde{E}_i)\geq s$, this
implies that ${\mathsf b}_{\mu_V}^{q}(\pi_V(E))\geq s$. Therefore,
we obtain ${\mathsf b}_{\mu_V}^{q}(\pi_V(E))\geq
\underline{\mathfrak{T}}_{\mu}^{q,m}(E)$.

\item First we suppose that
$\underline{\mathfrak{T}}_{\mu}^{q,m}(E)\geq m(1-q)$. If
$s>\underline{\mathfrak{T}}_{\mu}^{q,m}(E)$ we may cover $E$ by a
countable collection of sets $E_i$, which we may take to be compact,
such that $\underline{\tau}_{\mu}^{q,m}(E_i)<s$. By using Theorem
\ref{TH222} (3.) we have
$\underline{\tau}_{\mu_V}^{q}(\pi_V(E_i))\leq s$ for
$\gamma_{n,m}$-almost every  $V\in G_{n,m}.$ Proposition \ref{ah}
implies that ${\mathsf b}_{\mu_V}^{q}(\pi_V(E))\leq s$ for
$\gamma_{n,m}$-almost every  $V\in G_{n,m}$ and so, ${\mathsf
b}_{\mu_V}^{q}(\pi_V(E))\leq
\underline{\mathfrak{T}}_{\mu}^{q,m}(E)$ for $\gamma_{n,m}$-almost
every $V\in G_{n,m}.$ The proof of ${\mathsf
b}_{\mu_V}^{q}(\pi_V(E))\geq
\underline{\mathfrak{T}}_{\mu}^{q,m}(E)$ for all $V\in G_{n,m},$ is
identical to the proof of the above statement and is therefore
omitted.

Now, we suppose that $\underline{\mathfrak{T}}_{\mu}^{q,m}(E)\leq
m(1-q)$. It follows from Proposition \ref{ah} and Theorem
\ref{TH222} (3.) that ${\mathsf b}_{\mu_V}^{q}(\pi_V(E))=m(1-q)$
for $\gamma_{n,m}$-almost every  $V\in G_{n,m}.$
\end{enumerate}

\subsection{Proof of Theorem \ref{THbB}}
The proof of Theorem \ref{THbB} is straightforward from Remark
\ref{remark1} and \cite[Corollary 4.4]{O}.

\section{Application}\label{Sec6}
When  $\mu$ obeys the multifractal formalism over some interval, we
are interested in knowing whether or not this property is preserved
after orthogonal projections on $\gamma_{n,m}$-almost every linear
$m$-dimensional subspaces. In this section we  study the behavior of
projections of measures obeying to the multifractal formalism which
is based on the Hewitt-Stromberg measures. More specifically, we
prove that for $q>1$ if the multifractal formalism holds for $\mu$
at $\alpha=-{\mathsf B}'_{\mu}(q)$, it holds for $\mu_V$ for
$\gamma_{n,m}$-almost every $V\in G_{n,m}$. The multifractal
analysis is a natural framework to describe geometrically the
heterogeneity in the distribution at small scales of positive and
finite compactly supported Borel measures on $\mathbb{R}^n$.
Specifically, for such a measure $\mu$, this heterogeneity can be
classified by considering the iso-H\"{o}lder sets
$$
\mathscr{X}_{\mu}(\alpha)=\left\{x\in\supp\mu\;\Big|\;
\lim_{r\to0}\frac{\log\Big(\mu B(x,3r)\Big)}{\log r}=\alpha \right\}
$$
which form a partition of $\supp\mu$. Then the singularity spectrum
of the measure $\mu$ is the mapping $\alpha\mapsto \dim_{(.)}
\mathscr{X}_{\mu}(\alpha)$ and  this spectrum provides a geometric
hierarchy between the sets $\mathscr{X}_{\mu}(\alpha)$. We mention
that in the last decade there has been a great interest for the
multifractal analysis and positive results have been written in
various situations, see for example \cite{NB2, NB1, BBJ, BBH, OL,
OL1, O, SB3}. The function ${\mathsf B}_\mu(q)$ is related to the
multifractal spectrum of the measure $\mu$. More precisely,
$f^*(\alpha)=\inf_\beta\big(\alpha\beta+f(\beta)\big)$ denotes the
Legendre transform of the function $f,$ it has been proved in
\cite{NB1} a lower and upper bound estimate of the singularity
spectrum using the Legendre transform of the function ${\mathsf
B}_\mu(q)$. In this section, we will work with the following
formalism which a consequence of the multifractal formalism
developed in \cite{NB, NB1}.
\begin{theorem}\label{th3}
Let $\mu$ be a compactly supported Borel probability measure on
$\mathbb{R}^n$ and $ q\in\mathbb{R}$. Suppose that
\begin{enumerate}
\item there exists
a nontrivial (Frostman) measure $\nu_q$ satisfying
$$
\nu_q(B(x,r))\leq  \mu(B(x,3r))^q ~r^{{\mathsf
B}_\mu(q)}\quad\text{where}\quad x\in \supp\mu,\;\;0 < r < 1.
$$

\item ${\mathsf B}_{\mu}$ is differentiable at $q$.
\end{enumerate}
Then, for any $\alpha=-{\mathsf B}'_{\mu}(q)$
\begin{eqnarray*}
\dim_H
\mathscr{X}_{\mu}\big(\alpha\big)&=&\underline{\dim}_{MB}\mathscr{X}_{\mu}\big(\alpha\big)=\overline{\dim}_{MB}\mathscr{X}_{\mu}\big(\alpha\big)=\dim_P
\mathscr{X}_{\mu}\big(\alpha\big)\\& =&
B^*_{\mu}\big(\alpha\big)=b^*_{\mu}\big(\alpha\big)={\mathsf
b}^*_{\mu}\big(\alpha\big)={\mathsf B}^*_{\mu}\big(\alpha\big).
\end{eqnarray*}
\end{theorem}

The following proposition  has recently been obtained in \cite[Lemma
3.2]{O}.
\begin{proposition}\label{P}
Let $\mu$ be a compactly supported Borel probability measure on
$\mathbb{R}^n$. For $q\geq1$ and all $V\in G_{n,m}$, we have
 $$
b_{\mu_V}(q)\geq \max\Big(m(1-q), b_\mu(q)\Big).
 $$
\end{proposition}

 \bigskip
In the following, we study the validity of the multifractal
formalism under projection. More specifically, we obtain general
result for the multifractal analysis of the orthogonal projections
on $m$-dimensional linear subspaces of  measure $\mu$ satisfying the
multifractal formalism which is based on the Hewitt-Stromberg
measures.
\begin{theorem}
Let $\mu$ be a compactly supported Borel probability measure on
$\mathbb{R}^n$ and $ q>1$. Suppose that
\bigskip
\\$(\mathsf H_1)$ there exists
a nontrivial (Frostman) measure $\nu_q$ satisfying
$$
\nu_q(B(x,r))\leq  \mu(B(x,3r))^q ~r^{{\mathsf
B}_\mu(q)}\quad\text{where}\quad x\in \supp\mu,\;\;0 < r < 1,
$$
$(\mathsf H_2)$ ${\mathsf B}_{\mu}$ is differentiable at $q$,
 \bigskip
 \\$(\mathsf H_3)$ $({E}_i)_i$  be a cover
of $\supp\mu$ by a countable collection of compact sets is such that
$b_{\mu}^q(E_i\cap\supp\mu)\geq \max(-m, m(1-q))$ for all $i$.
 \bigskip
\\Then, for any $\alpha=-{\mathsf B}'_{\mu}(q)$ and  $\gamma_{n,m}$-almost every $V\in
G_{n,m}$,
\begin{eqnarray*}
&\underline{\dim}_{MB}\mathscr{X}_{\mu}\big(\alpha\big)=\overline{\dim}_{MB}\mathscr{X}_{\mu_V}(\alpha)=
\underline{\dim}_{MB}\mathscr{X}_{\mu_V}\big(\alpha\big)=\overline{\dim}_{MB}\mathscr{X}_{\mu}(\alpha)=\\&\dim_H
\mathscr{X}_{\mu}\big(\alpha\big)=\dim_P
\mathscr{X}_{\mu_V}\big(\alpha\big)= \dim_H
\mathscr{X}_{\mu_V}\big(\alpha\big) = \dim_P
\mathscr{X}_{\mu}\big(\alpha\big) =\\& B^*_{\mu}\big(\alpha\big) =
b^*_{\mu}\big(\alpha\big)={\mathsf
b}^*_{\mu}\big(\alpha\big)={\mathsf B}^*_{\mu}\big(\alpha\big).
\end{eqnarray*}
\end{theorem}
\noindent {\bf Proof.} It follows from Theorems \ref{TH2} and
\ref{TH3}, Propositions \ref{P} and \ref{L1}, $(\mathsf H_1)$ and
$(\mathsf H_3)$ that, for $\gamma_{n,m}$-almost every $V\in
G_{n,m}$,
\begin{equation}\label{qq}
B_{{\mu}_V}({q})={\mathsf b}_{{\mu}_V}({q})={\mathsf
B}_{{\mu}_V}({q})=b_{{\mu}_V}({q})= b_{{\mu}}({{q}})={\mathsf
b}_{{\mu}}({{q}})={\mathsf B}_{{\mu}}({{q}})=B_{{\mu}}({q}).
\end{equation}
The hypothesis $(\mathsf H_1)$,  (\ref{qq}) and the proof of Lemma
3.2 in \cite{O} ensure that, there exists a positive constant  $c$
such that
\begin{center}
$0<\mathscr{H}_{\mu}^{q,{\mathsf B}_{\mu}(q)}(\supp\mu)\leq
c~\mathscr{H}_{\mu_V}^{q,{\mathsf B}_{\mu_V}(q)}(\supp\mu_V),$\quad
for $\gamma_{n,m}$-almost every $V\in G_{n,m}$.
\end{center}
By using \cite[Theorem 5.1]{O} there exists a Frostman measure
$\nu_q$ satisfying
$$
\nu_q(B(y,r))\leq  \mu_V(B(y,3r))^q ~r^{{\mathsf
B}_{\mu_V}(q)}\quad\text{where}\quad y\in \supp\mu_V,\;\;0 < r < 1.
$$
Now, it follows from the hypothesis $(\mathsf H_2)$  and Theorem
\ref{th3} that
\begin{eqnarray}\label{4}
\underline{\dim}_{MB}\mathscr{X}_{\mu_V}(\alpha)\geq\dim_H\mathscr{X}_{\mu_V}\big(\alpha\big)\geq
q\alpha+{\mathsf B}_{\mu}(q),
\end{eqnarray}
$\text{for}\;\gamma_{n,m}\text{-almost every} \;V\in G_{n,m}$.
Clearly the assumption (\ref{qq}) give that
\begin{eqnarray}\label{3}
\overline{\dim}_{MB}\mathscr{X}_{\mu_V}(\alpha)=\dim_P\mathscr{X}_{\mu_V}\big(\alpha\big)&\leq&
{\mathsf B}^*_{\mu_V}\big(\alpha\big)= {\mathsf
B}^*_{\mu}\big(\alpha\big),
\end{eqnarray}
for $\gamma_{n,m}$-almost every $V\in G_{n,m}$. Now combining
\eqref{4}, \eqref{3} and  the equalities \eqref{qq} we get the
desired result.

\bigskip\bigskip
\noindent {\bf Bilel SELMI}\\ \\
\noindent Analysis, Probability and Fractals Laboratory LR18ES17\\
Faculty of Sciences of Monastir\\
Department of Mathematics\\
University of Monastir\\
5000-Monastir\\
Tunisia\\ \\
E-mail:\; bilel.selmi@fsm.rnu.tn
\end{document}